\documentclass[twocolumn,amsmath,amssymb,aps,secnumarabic,%
    nofootinbib,groupedaddress,floatfix]{revtex4}

\usepackage{mathpazo,graphics,amssymb,amsmath,amsthm}
\usepackage{times}

\newtheorem{theorem}{Theorem} 

\newtheorem{lemma}[theorem]{Lemma}

\newtheorem{corollary}[theorem]{Corollary} 

\newcommand{\R}{\mathbb{R}}

\newcommand{\Z}{\mathbb{Z}}
\renewcommand{\d}{\partial}

\newcommand{\bq}{\begin{equation}}
\newcommand{\eq}{\end{equation}}

\newcommand{\abs}{\mathrm{abs}}
\newcommand{\Ie}{\emph{I.e.}}
\newcommand{\ie}{\emph{i.e.}}
\newcommand{\grad}{\vec{\nabla}}
\renewcommand{\div}{\vec{\nabla} \cdot}
\newcommand{\vx}{\vec{x}}

\newcommand{\fullfigure}[2]{
    \begin{figure}[htb] \begin{center} \includegraphics{#1.eps}
    \end{center} \caption{#2}\label{#1} \end{figure}}
    
\newcommand{\widefullfigure}[2]{
    \begin{figure*}[tb] \begin{center} \includegraphics{#1.eps}
    \end{center} \caption{#2}\label{#1} \end{figure*}}

\begin{document}
\author{Greg Kuperberg}
\email{greg@math.yale.edu}
\affiliation{Department of Mathematics, Yale University, New Haven, CT 60637}
\thanks{The author was supported by an NSF
    Postdoctoral Fellowship, grant \#DMS-9107908.}
\title{A volume-preserving counterexample to the Seifert conjecture}
\date{April 29, 1995}

\begin{abstract}
We prove that every 3-manifold possesses a $C^1$, volume-preserving
flow with no fixed points and no closed trajectories.  The main construction is a
volume-preserving version of the Schweitzer plug.  We also prove that
every 3-manifold possesses a volume-preserving, $C^\infty$ flow with
discrete closed trajectories and no fixed points (as well as a PL flow with the
same geometry), which is needed for the first result.  The proof uses a
Dehn-twisted Wilson-type plug which also preserves volume.

\vspace{\baselineskip} \noindent MSC classification (1991): Primary 58F25;
Secondary 57R30, 28D10, 58F05, 58F22
\end{abstract}
\maketitle

\begin{theorem} Every 3-manifold possesses a $C^1$ volume-preserving
flow with no fixed points and no closed trajectories. \label{thmain}
\end{theorem}

The author was motivated to consider Theorem~\ref{thmain} by the recent
discovery of a real analytic counterexample to the Seifert conjecture
\cite{kseifert}; the conjecture states that every flow on $S^3$ has either a
fixed point or a closed trajectory.  However, the construction
presented here is
based on the original $C^1$ counterexample due to Schweitzer
\cite{Schweitzer} and not the new counterexample.

An important property of a volume-preserving flow in 3 dimensions with no fixed
points is that the parallel 1-dimensional foliation is transversely symplectic. 
In particular, such a flow on a 3-manifold $M$ can be understood as a
Hamiltonian flow coming from a symplectic structure on $M \times \R$.  In this
context, the Weinstein conjecture \cite{Weinstein} provides an interesting
contrast to Theorem~\ref{thmain}.  It states that a flow on a closed
$(2n+1)$-manifold $M$ which is not only transversely symplectic but also contact
must have a closed trajectory if $H^1(M,\Z) = 0$.  (A contact form on an
$(2n+1)$-manifold
is a 1-form $\omega$ such that $\omega \wedge (d \omega)^{\wedge n}$ does not
vanish; the corresponding contact flow is parallel to the kernel of $d\omega$.)
In particular, Hofer \cite{Hofer} has recently established the Weinstein
conjecture for $S^3$ and his result applies to the $C^1$ category.  Thus, the
flow established by theorem~\ref{thmain} is not contact.

For most manifolds, although not $S^3$, Theorem~\ref{thmain} depends on the
following result.

\begin{theorem} Every 3-manifold possesses a $C^\infty$ volume-preserving flow
with no fixed points and a discrete collection of closed trajectories, as
well as a transversely measured 1-dimensional PL foliation with discrete
closed leaves.
\label{thdiscrete} \end{theorem}

Theorem~\ref{thdiscrete} is an extension of the 3-dimensional case of Wilson's
theorem~\cite{Wilson}, which establishes flows with no fixed points and discrete
closed trajectories, but without the volume preservation condition.  Like 
Wilson's theorem
and all known counterexamples to the Seifert conjecture, Theorems~\ref{thmain}
and~\ref{thdiscrete} both use the standard technique of constructing plugs and
inserting them into other flows.  However, volume preservation restricts the
behavior of a plug, and in particular a volume-preserving plug cannot stop an
open set.  To work around this serious constraint, Theorem~\ref{thdiscrete}
uses twisted plugs, which are plugs whose insertions into manifolds
can change the topology of the manifolds.  However, the
twisted plugs constructed here are only $C^\infty$ and not real
analytic.  The generalization of Theorem~\ref{thdiscrete} to the
real analytic case remains open.

I am grateful to Krystyna Kuperberg for encouragement and important discussions
about the research presented here.  I would also like to thank Shmuel Weinberger
and \'Etienne Ghys for their interest in the results and useful comments.

\section{\label{sintro} Preliminaries}

In this paper, we will mainly consider three smoothness categories:
$C^r$ for finite $r$, $C^\infty$ or smooth, and PL.  In many contexts,
an object will have implicit smoothness; for example, a map between
PL manifolds will be assumed to be PL.  Unless explicitly stated otherwise,
all of the arguments assume that manifolds are oriented and connected but
generalize easily to non-orientable and disconnected manifolds.

The paper will use several different $C^\infty$ bump functions and
transition functions.  Let $b:[0,1] \to \R$ be a non-negative $C^\infty$
function with support $[1/3,2/3]$ whose integral is 1 and which does
not exceed 4.  Let $B:[0,1] \to \R$ be
a non-negative $C^\infty$ function whose value and derivatives vanish at $0$ and
$1$ and such that $B(x) > b(x) \ge 0$ for $ 0 < x <1$.  Let $e:[-1,1] \to \R$ be
a non-negative $C^\infty$ function such that:
$$ e(x) \left\{\begin{array}{lr}
= 0 & |x| \ge 2/3 \\
< 1 & |x| > 1/3 \\
= 1 & |x| \le 1/3
\end{array}\right. $$
Finally, let $o:[-1,1] \to \R$ be a $C^\infty$ odd, increasing function (\ie,
$o(x) > o(y)$ if $x>y$ and $o(-x) = -o(x)$) such that $o(1) = 1$ and all
derivatives of $o$ vanish at the origin.

\subsection{Foliations and volume-preserving flows}

In a standard mathematical treatment of fluid motion, a vector field $\vec{v}$
in $\R^n$ represents a static flow, and if the divergence equation
$$\div \vec{v} = 0$$
holds, the flow preserves volume.  This definition must be carefully
generalized to flows on manifolds. A {\em smooth measure structure} on a
smooth $n$-manifold is in general given by a smooth, non-vanishing $n$-form
or {\em volume form}. It is a simple result that any volume form is
equivalent to Lebesgue measure by a local diffeomorphism.  More
interestingly, Moser's theorem~\cite{Moser:1965} states that a compact
manifold $M$ with two volume forms $\mu_1$ and $\mu_2$ with the same total
volume admits a diffeomorphism taking $\mu_1$ to $\mu_2$.  Given a volume
form $\mu$ on a manifold $M$, the divergence equation for a tangent vector
field $\vec{v}$ on a manifold becomes
$$d(\iota_{\vec{v}}(\mu)) = 0, \label{diveq}$$
where the operator $\iota_{\vec{v}}$ is contraction with $\vec{v}$.  The closed
$(n-1)$-form $\iota_{\vec{v}}(\mu)$ is a {\em flux form}.  This formalism is
compatible with another view of flows. For $M$ closed, a flow can be defined as
a smooth group action $\Phi: \R \times M \to M$; the vector field $\vec{v}$ is
related to the group action $\Phi$ by
$$\vec{v} = \frac{d \Phi}{d t}.$$
It is easy to check that the condition that $\mu$ is invariant under
$\Phi$ is equivalent to the divergence equation.

Since contraction with $\mu$ is an invertible linear transformation, given a
flux form $\omega$, any volume form $\mu$ yields a vector field $\vec{v}$
such that 
$$\omega = \iota_{\vec{v}}(\mu),$$ 
with the conclusion that $\vec{v}$ preserves $\mu$. Moreover, the trajectories
of $\vec{v}$, that is the curves parallel to $\vec{v}$, are determined by
$\omega$, since $\omega$ has a 1-dimensional kernel at a point of $M$ where
it does not vanish, and it is easy to check that $\vec{v}$ is a non-vanishing
vector in the kernel at that point.  Where $\omega$ does vanish, $\vec{v}$
vanishes also.

A useful special case of the flux form view is $n=2$, for then a flux form
$\omega$ is the differential of a possibly multiply-valued function $f$. 
Dualizing by some area form, the trajectories of the vector field $\vec{v}$
obtained from $\omega = df$ are simply the contours of $f$.  If the manifold
$M$ has a Riemannian metric and the volume form is given by the metric, then
$\vec{v}$ can also be defined as $J(\grad f)$, where $J$ is a rotation by 90
degrees. This expression is also validated by the standard identity 
$$\div J(\grad f) = 0$$
from 2-dimensional vector calculus.

If $\vec{v}$ has no fixed points, its trajectories form a 1-dimensional,
oriented foliation $\cal F$ of the manifold $M$.  The usual Seifert conjecture
can be phrased as a question about foliations rather than flows: Is there an
oriented 1-foliation of $S^3$ with no closed leaves? In the volume-preserving
case, the foliation $\cal F$ is determined by the flux form $\omega$ and
otherwise does not depend on $\vec{v}$ or the volume form; $\cal F$ can be
defined as the unique foliation parallel to the kernel of $\omega$, which is a
line bundle over $M$.  Since $\omega$ is closed, it determines a (transverse)
measure on $\cal F$.  A {\em measure} on a $k$-foliation of an $n$-manifold is
in general defined as a measure for every transverse $(n-k)$-disk which is
invariant under isotopy of the disk parallel to the foliation. In this case,
given an $(n-1)$-disk $D^{n-1}$ transversely embedded by $\alpha:D^{n-1} \to M$,
the measure on $D^{n-1}$ is given by the pullback $\alpha^*(\omega)$, which is a
smooth volume form.  The measure induced by $\omega$ can be called smooth and
locally Lebesgue just as volume forms are.

As a point of terminology, if $\cal F$ is a foliation of $M$, $M$
is the {\em support} of $\cal F$.  Also, henceforth the term foliation will
mean an oriented 1-foliation except where explicitly stated otherwise.

The simplest important construction with foliations is the suspension. (In
topology the suspension is called the mapping torus and the term suspension
is used for a different construction). Given a manifold $M$ and a
diffeomorphism $\sigma: M \to M$, the {\em suspension} is the manifold $M \times
\R$ foliated by lines $p \times I$ and quotiented by the relation $(p,x)
\sim (\sigma(p),x-1)$. The main properties of the suspension used in this
paper are that $\cal F$ is measured if $\sigma$ preverses volume of $M$ and
that closed leaves of $\cal F$ correspond to finite orbits of $\sigma$.  Note
also that if $\sigma$ is isotopic or even pseudoisotopic to the identity,
the support of $\cal F$ is diffeomorphic to $M \times S^1$.

\subsection{PL foliations}

Recall that a $k$-foliation of an $n$-manifold is an atlas of charts
such that each gluing map preserves horizontal $k$-planes in $\R^n$.
In other words, each gluing map can be written as
\begin{multline*}
g(\vx) = (g_1(\vx),g_2(\vx),\ldots,g_k(\vx), \\
g_{k+1}(x_{k+1},\ldots,x_n),\ldots,g_n(x_{k+1},\ldots,x_n)).
\end{multline*}
A $k$-foliation is PL if the gluing maps are PL, smooth if the gluing maps are
smooth, and has a (locally Lebesgue) measure if the transverse part
$(g_{k+1},\ldots,g_n)$ of each gluing map preserves Lebesgue measure on
$\R^{n-k}$.  In the smooth case, this definition is equivalent to the one in
terms of flux forms. However, in the PL case, the atlas definition seems to be
the best substitute for flux forms.  Note that suspensions generalize the PL
case.

\widefullfigure{}{Transfer of measure between adjacent simplices}

Similarly, a (locally Lebesgue) measure structure on an $n$-manifold
can be defined as an atlas of charts such that the gluing maps preserve
Lebesgue measure on $\R^n$; in the smooth case such an atlas is
equivalent to a volume form.  In the PL case, such an atlas is an 
example of a simplicial measure.  A measure on a PL manifold is
{\em simplicial} relative to a triangulation $T$ if on each simplex
the measure is given by a linear embedding of the simplex in Euclidean
space.  The following analogue of Moser's theorem demonstrates
that simplicial measures are the PL analogue of volume forms:

\begin{theorem} Two simplicial measures on a connected, compact PL
$n$-manifold $M$ with the same total volume are equivalent by a PL
homeomorphism.  Moreover, any simplicial measure is locally PL-Lebesgue.
\label{thplmoser}
\end{theorem}
\begin{proof} Let $\mu_1$ and $\mu_2$ be two simplicial measures on $M$ and
let $\cal T$ be a triangulation for which both measures are simplicial.  For
any pair of simplices $T_1$ and $T_2$ of $\cal T$ that meet at an
$n-1$-dimensional face, there is a family of PL homeomorphisms of $T_1 \cup
T_2$ which take measures which are simplicial on the triangulation
$\{T_1,T_2\}$ to other such measures and which transfer measure from $T_1$ to
$T_2$.  Figure~\ref{ftransfer} shows an example of such a homeomorphism: In
the figure, $T_1$ and $T_2$ are embedded in such a way that their volumes are
proportional to their given measure.  We retriangulate $T_1 \cup T_2$ with
simplices $U_1,\ldots,U_n$ that share a 1-dimensional edge.  There is then a
homeomorphism $\phi$ which is linear and volume-preserving on each simplex
$U_i$ such that image is a union of two simplices $T'_1$ and $T'_2$ whose
volumes differ from $T_1$ and $T_2$. By identifying $T'_1$ with $T_1$ and
$T'_2$ with $T_2$, $\phi$ can be understood as a PL map that transfers
measure from $T_1$ to $T_2$.

Thus, using PL homemorphisms modelled on $\phi$, we can transfer
measure arbitrarily between adjacent simplices of $\cal T$, as
long as the measure of each simplex remains positive.  Such
moves clearly suffice to connect any two measures $\mu_1$ and $\mu_2$:
By analogy, a connected graph of people, each with a positive
amount of money, can arbitrarily redistribute their assets
solely by having adjacent individuals transfer money; moreover,
the transfers need not drive any individual into debt.

The PL maps so produced fix the vertices of $\cal T$, even though they are
not linear on the simplices of $\cal T$.  To show the second claim, let $\mu$
be a measure simplicial relative to $\cal T$, let $x$ be a point in $M$, and
let $y$ be a point in $M$ in the interior of some simplex of $\cal T$. The
measure $\mu$ is clearly PL-Lebesgue in a neighborhood of $y$.  Let $\mu_1 =
\mu$, let $\mu_2 = \alpha(\mu_1)$, where $\alpha$ is a PL homeomorphism of
$M$ that takes $y$ to $x$, and let $\cal U$ be a mutual refinement of $\cal
T$ and $\alpha({\cal T})$.  The points $x$ and $y$ are both vertices of $\cal
U$ and $\mu_1$ and $\mu_2$ are both simplicial relative to $\cal U$. 
Applying the above argument, $\mu$ is locally the same at $x$ and $y$, and
therefore $\mu$ is PL-Lebesgue at $x$ also.
\end{proof}

\subsection{The $C^r$ case}
\label{scrcase}

A non-vanishing $C^r$ vector field on a manifold $M$ yields a $C^r$ foliation
(a foliation with $C^r$ gluing maps), which yields a $C^r$ structure for the
manifold, but unfortunately the smoothness of vector fields on a $C^r$
manifold is only defined up to $C^{r-1}$.  Similarly, a $C^r$ manifold with
volume-preserving gluing maps only has a $C^{r-1}$ volume form.  If both
structures are present, the vector field can be smoothed to $C^r$
\cite{Thurston}, but usually at the expense of crumpling a $C^\infty$ volume
form to $C^{r-1}$.  Alternatively, by a refinement of Moser's theorem
\cite{Moser:1981}, a $C^{r-1}$ plus H\"older volume form can be smoothed to
$C^\infty$, but only by a $C^r$ plus H\"older diffeomorphism which might
crumple a $C^r$ vector field so that it is only $C^{r-1}$ plus H\"older.

In this paper, a volume-preserving $C^r$ flow means a $C^r$ vector
field which is divergenceless relative to a smooth volume form on a smooth
manifold.  It suffices to consider $C^r$ flux forms on smooth manifolds (or
at least $C^{r+1}$ manifolds). In particular, such a flux form defines a
measured $C^r$ foliation.  Although a $C^r$ flux form is not exactly the same
as such a foliation, it is very similar and for some purposes it will be
convenient to treat it as one.  When we need to glue together flux forms, we
will require that they are smooth in the gluing regions, so that in these
regions they are equivalent to smooth measured foliations.

\fullfigure{}{Corner separation}

\section{Plugs}
\label{splugs}

To define plugs, we must consider a class of manifolds with at least some
kinds of corners.  The smallest convenient such class is the class of orthant
manifolds.  An $n$-dimensional {\em orthant manifold} is a Hausdorff space
locally homeomorphic to some open subset of the orthant in $\R^n$ of points
with non-negative coordinates.  In the PL category, an orthant manifold is
just a manifold with boundary, but in the smooth category, the boundary might
not be smooth. For example, a parallelopiped is a smooth orthant manifold.

A foliation can have many different kinds of structure at the boundary of an
orthant manifold, not to mention the boundary of an ordinary manifold with
boundary, but to define a plug three kinds of boundary structure suffice:  {\em
parallel boundary}, {\em transverse boundary}, and  {\em corner separation}
between parallel and transverse boundary.  Figure~\ref{fcorner} shows an example
of each type of boundary. Recall some definitions from reference \cite{kk}: A
{\em flow bordism} is a foliation $\cal P$ on a compact orthant manifold $P$
such that $\d P$ is entirely transverse boundary, parallel boundary, or corner
separation, and such that all leaves in the parallel boundary of $P$ are
finite.  If $\cal P$ is a flow bordism, let $F_-$ be the (closure of) all
transverse boundary oriented inward, and similarly let $F_+$ be the transverse
boundary oriented outward.  The foliation $\cal P$ might in addition have one or
both of the following properties:
\begin{itemize}
\item[(i)] There exists an infinite leaf with an endpoint in $F_-$.
\item[(ii)] There exists a manifold $F$ and homeomorphisms
$\alpha_\pm:F \to F_\pm$ such that if $\alpha_+(p)$ and $\alpha_-(q)$
are endpoints of a leaf of $\cal P$, then $p=q$.
\end{itemize}
If $\cal P$ satisfies property (ii), it has {\em matched ends}. The foliation
$\cal P$ is a {\em plug} if it has properties (i) and (ii), but only a {\em
semi-plug} if it has property (i) but not property (ii).  It is an {\em
un-plug} if has property (ii) but not property (i).  The manifold $F_-$ is
the {\em entry region} of $\cal P$, while $F_+$ is the {\em exit region}. If
$\cal P$ has matched ends, then $F$ is the base of $\cal P$.   The {\em entry
stopped set} $S_-$ of $\cal P$ is the set of points of $F_-$ which  are
endpoints of infinite leaves; the exit stopped set $S_+$ is defined
similarly.  If $\cal P$ has matched ends, the stopped set $S$ is defined as
$\alpha_-^{-1}(S_-) = \alpha_+^{-1}(S_+)$.  If $\cal P$ has matched ends and
$S$ contains an open set, then $\cal P$ {\em stops content}, \ie, has
wandering points in $F$.

An important construction due to Wilson \cite{Wilson} turns a semi-plug into a
plug.  If ${\cal P}_1$ and ${\cal P}_2$ are two flow bordisms such that the exit
region of ${\cal P}_1$ is the same as the entry region of ${\cal P}_2$, their
{\em concatenation} is a flow bordism obtained by identifying trivially foliated
neighborhoods of this shared region.  The {\em mirror image} $\bar{\cal P}$ of a
flow bordism $\cal P$ is given by reversing the orientation of the leaves of
$\cal P$, which has the effect of switching the entry and exit regions.  The
{\em mirror-image construction} is the concatenation of $\cal P$ and $\bar{\cal
P}$; it is easy to see that the result of this concatenation has matched ends.

The primary purpose of plugs is the operation of insertion.  An {\em insertion
map} for a plug $\cal P$ into a foliation $\cal X$ is an embedding $F \to X$ of
the base of $\cal P$ which is transverse to $\cal X$.  Such an insertion map can
be extended to an embedding $\sigma:F \times I \to X$ which takes the fiber
foliation of $F \times I$ to $\cal X$.  An $n$-dimensional plug $\cal P$ is
insertible if $F$ admits an embedding in $\R^n$ which is transverse to vertical
lines.  Such an embedding is equivalent to a bridge immersion of $F$ in
$\R^{n-1}$, \ie, an immersion which lifts to an embedding one dimension higher. 
Figure~\ref{fbridge} shows a bridge immersion of a punctured torus $pT$; the
corresponding embedding of $F \times I$ is the one that Schweitzer also uses.

\fullfigure{}{A bridge immersion}

Let $N_{F \times I}$ be an open neighborhood of $\d(F \times I)$. The next step
in plug insertion is to remove $\sigma((F \times I)-N_{F \times I})$ from $X$
and glue the open lip $\sigma(N_{F \times I})$ to $P$ by a leaf-preserving
homeomorphism $\alpha:N_{F \times I} \to N_P$, where $N_P$ is a neighborhood of
$\d P$.  Moreover, the identification $\alpha$ should satisfy $\alpha(p,0) =
\alpha_-(p)$ and $\alpha(p,1) = \alpha_+(p)$.  A map $\alpha$ with these
properties is an {\em attaching map} for $\cal P$.  As explained in Reference~\cite{kk},
plugs always possess attaching maps.

Let $\sigma$ be an insertion map of a plug $\cal P$ into a foliation $\cal X$ on
a manifold $X$, and let $\hat{\cal X}$ be the foliation on the manifold
$\hat{X}$ resulting from the insertion of $\cal P$ into $\cal X$.  The plug
$\cal P$ is {\em untwisted} if the attaching map $\alpha$ extends to a
homeomorphism $F \times I \to P$, and {\em twisted} otherwise.  If $\cal P$ is
untwisted, then $X$ and $\hat{X}$ are necessarily homeomorphic, while if $\cal
P$ is twisted, then $X$ and $\hat{X}$ need not be homeomorphic.  This paper will
use both twisted and untwisted plugs.

A useful lemma about plugs proved in Reference~\cite{kk} is the
following:

\begin{lemma} A flow bordism with an infinite leaf with non-empty entry or exit
region is either a plug or a semi-plug. \label{linfleaf}
\end{lemma}

\subsection{Measured and $C^r$ plugs}
\label{scrplug}

The technique of plugs generalizes without any substantive changes to the
category of measured foliations, either smooth or PL.  The base $F$ of a
measured plug is measured, but by Moser's theorem, the only relevance of this
structure is that a bridge immersion of $F$ with large volume into a disk
with small volume is inconvenient, although not strictly impossible, if $F$
has large volume.  One way to overcome this inconvenience is to rescale the
transverse measure of the plug to make the measure of $F$ small.  Note
also that a measured plug cannot stop content.

The category of measured, $C^r$ foliations is trickier.  Following the
prescription of subsection~\ref{scrcase}, a measured $C^r$ flow bordism is
realized by a $C^r$ flux form on a smooth manifold.  A flow bordism with support
$P$ is {\em attachable} if the flux form is smooth in a neighborhood $N_P$ of the
boundary, so that the foliation method can be used to insert it without loss of
smoothness.

A $n$-dimensional, measured $C^r$ semi-plug $\cal P$ with support $P$ can
always be made attachable by the following method:  Since the flux form
$\omega$ is defined over all of $P$ and $P$ bounds $\d P$, it is the
differential of an $(n-2)$-form $\nu$ in a neighborhood of the boundary
$N_P$.  Let $\nu'$ be an $(n-2)$-form which is a smooth approximation to
$\nu$ in a smaller neighborhood of $\d P$, agrees with $\nu$ in a
neighborhood of $P - N_P$, and is an interpolation with a smooth bump
function in between. In addition, choose $\nu'$ so that $d\nu'$ has the same
parallel and transverse boundary at $\d P$ as does $\omega$.  Then the flux
form which is $d\nu'$ on $N_P$ and $\omega$ on $P - N_P$ yields an attachable
semi-plug $\cal P'$ with the same leaf structure as $\cal P$. Furthermore,
the mirror-image construction applied to $\cal P'$ yields an attachable plug.

\section{\label{sdiscrete} Isolated closed trajectories}

The main construction of the proof of Theorem~\ref{thdiscrete} is a measured,
Dehn-twisted plug.  Before constructing or even defining such a plug, we recall
several facts about Dehn twists and Dehn surgery:  The boundary of a solid torus
$S^1 \times D^2$ has a distinguished embedded circle, the {\em meridian},
which is unique
up to isotopy and which is identified by the fact that it bounds a disk in the
solid torus.  A {\em framing} of a solid torus is a homotopy class of another
circle in the boundary; the framing may or may not equal the meridian.  A
framing is {\em integral} if it homologically crosses the meridian exactly once.
A {\em Dehn surgery} on a 3-manifold consists of removing a collection of disjoint
framed tori and gluing them back in such a way that the new meridian circles
match the old framing circles; the topology of the resulting manifold does not
otherwise depend on the gluing maps.  The Lickorish-Wallace theorem
\cite{Rolfsen} asserts that every closed, oriented 3-manifold can be obtained
from $S^3$ by integral Dehn surgery, or equivalently every closed, oriented
3-manifold can be obtained from every other by integral Dehn surgery. In the
second formulation, the theorem also holds for non-orientable manifolds.

Suppose that $\cal D$ is a plug with base $F = S^1 \times I$ whose support $P$
is homeomorphic to a solid torus $S^1 \times D^2$. Recall that there is an
attaching map $\alpha:N_{F \times I} \to N_P$ between neighborhoods of the
boundary, and note that the thickened base $F \times I$ is also a solid torus. 
If $p \in S^1$, the curve $m = \{p\} \times \d(I \times I)$ is a meridian of $F
\times I$, while the curve $l = S^1 \times \{0\} \times \{0\}$ is a convenient
standard framing which might be called the {\em longitude}. Recall that when $\cal D$
is inserted, its support $P$ replaces an image of $F \times I$ by the attaching
map $\alpha$.  Therefore if $\alpha(m)$ is a meridian of $P$, meridian replaces
meridian, $\alpha$ extends to a homeomorphism $\alpha:F \times I \to P$, and
$\cal D$ is untwisted. If, alternatively, the meridian of $P$ replaces some
other curve of $F \times I$, $\cal D$ can be called Dehn-twisted, because its
insertion effects a Dehn surgery. In particular, if either $\alpha(m+l)$ or
$\alpha(m-l)$ (using homological notation for other curves besides $m$ and $l$)
is a meridian of $P$, $\cal D$ is integrally Dehn-twisted.

An integrally Dehn-twisted plug $\cal D$, assuming that it exists, can be used
to construct a foliation on any closed, oriented 3-manifold with finitely
many closed leaves as follows:  The 3-torus $T^3$ possesses a smooth, measured
foliation $\cal T$ such that all leaves are dense: If $T^3$ is given with
periodic coordinates $\theta_1, \theta_2, \theta_3$, define $\cal T$ to be
parallel to the vector field
$$r_1 \frac{\d}{\d\theta_1} + r_2 \frac{\d}{\d\theta_2}
+ r_3\frac{\d}{\d\theta_3},$$
where $r_1$, $r_2$, and $r_3$ are linearly independent over the rationals. Let
$M$ be some other 3-manifold, and let $L$ be a link in $T^3$ such that some
integral surgery on $L$ yields $M$.  If the link $L$ is transverse to $\cal T$,
which can always be achieved by isotopy, then $L$ can be extended to an
insertion map for several copies of $\cal D$. The longitudes of the thickened
bases $F \times I$ along $L$ are determined by $\cal T$; they can be chosen to
be any desired integral framing by adding coils to $L$, as shown in
Figure~\ref{fcoils}.  The framing for the surgery induced by inserting $\cal D$
is then given by the formula $m\pm l$ above, and this is also an arbitrary
integral framing on each component of $L$.

\fullfigure{}{Coiling a Dehn-twisted insertion}

Non-orientable manifolds can similarly be handled as follows: A rotation of a
round 2-sphere $S^2$ by an irrational angle descends to a volume-preserving,
smooth diffeomorphism of the projective plane $\R P^2$ with only one periodic
point, a fixed point.  The suspension of this diffeomorphism is therefore a
measured
foliation of $\R P^2 \times S^1$ with one closed leaf.  Every other non-orientable,
closed 3-manifold can be obtained from this one by appropriate insertions of
$\cal D$.  An alternative approach is to use the wormhole plug defined in
subsection~\ref{snoncompact} to add a non-orientable handle to a foliated,
orientable 3-manifold.

In conclusion, the smooth, compact case of Theorem~\ref{thdiscrete} follows
from the following lemma:

\begin{lemma} 
There exists a smooth, measured, integrally Dehn-twisted plug $\cal D$
with two closed leaves. \label{ldehn}
\end{lemma}
\begin{proof} As a warm-up, we construct an untwisted, measured plug with two
closed leaves.  Let $F = \{(r,\theta) | 1 \le r \le 3\}$ be an annulus in the
plane given in polar coordinates, but with the volume form $dr \wedge d\theta$
rather than the form given by the embedding in the plane.  Consider $C = F
\times [-1,1]$ in cylindrical coordinates $r$, $\theta$, and $z$.  Let $f:[1,3]
\times [-1,1] \to \R$ be given by
$$f(r,z) = z^2(r-2) + (1-z^2)(r-2)^3.$$
The contours of $f$ are given in Figure~\ref{ffcontours}.  The function $f$ has
one critical point at $(2,0)$, and all contours of $f$ connect the top and the
bottom, although the $r=2$ contour is singular. Let $\vec{W}$ be a vector field
on $C$ given by
$$\vec{W_s} = J(\grad f) + \frac{\d}{\d \theta}.$$
Let ${\cal W}_s$ be the foliation of $C$ which is parallel to $\vec{W}$. The
vector field $\vec{W_s}$ is divergenceless because both terms are
divergenceless, and therefore ${\cal W}_s$ is measured. By the geometry of $f$,
the leaves of ${\cal W}_s$ at $r \ne 2$ connect the top and the bottom of $C$,
but the leaves at $r = 2$ spiral to a closed leaf with $r=2$ and $z = 0$. It
is easy to check that $(\d F) \times I$ is parallel boundary of ${\cal W}_s$,
while $F \times \d I$ is transverse boundary. In conclusion, ${\cal W}_s$ is a
semi-plug with one closed leaf. The mirror-image construction described in
section~\ref{sintro} applied to ${\cal W}_s$ yields a plug $\cal W$ with two
closed leaves.

\fullfigure{}{Contours of $f$}

The plug $\cal W$ is necessarily untwisted, because in the notation preceding
the lemma, the circle $\alpha(c)$ consists of two arcs with constant $\theta$,
one in the entry region and the other in the exit region, connected by two
leaves of $\cal W$.  (See Figure~\ref{fmeridian}a.)  By the mirror-image
construction, if either leaf winds by some angle $\theta$ in ${\cal W}_s$, it
unwinds by the same angle in the mirror image $\bar{{\cal W}_s}$, so the two
leaves together with the two connecting arcs do not wind around the hole of the
support of $\cal W$ and $\alpha(c)$ is a meridian.

\widefullfigure{}{The curves $\alpha(c)$ in $\cal W$ and $\cal P$}

The main construction is a variant $\cal P$ which is a concatenation of a
semi-plug ${\cal P}_1$ and the mirror image of another semi-plug ${\cal P}_2$
which is a modification of ${\cal P}_1$.  Both semi-plugs are supported on the
space $C$ defined above.  The semi-plug ${\cal P}_2$ is parallel to the vector
field
$$\vec{P_2} = J(\vec{\nabla}g) + \frac{\d}{\d \theta},$$
where $g:[1,3] \times [-1,1] \to \R$ is given by
$$g(r,z) = e(z)o(r-2) + (1-e(z))(r-2).$$
The function $g$ has zero first derivative on the line segment
$\{2\} \times [-1/3,1/3]$, and therefore the foliation ${\cal P}_2$
has an annulus of closed leaves $\{2\} \times S^1 \times [-1/3,1/3]$.
The semi-plug ${\cal P}_1$ is parallel to the vector field 
$$\vec{P_1} = J(\vec{\nabla}g) +
(3\pi b(\frac{1+3z}{2})o'(r-2)+1)\frac{\d}{\d \theta}$$
when $z \in [-1/3,1/3]$ and $r > 2$, and equals $\vec{P_2}$ otherwise. The
coefficient of $\frac{\d}{\d \theta}$, although complicated, does not involve
$\theta$, so $\vec{P_1}$ is still a sum of two divergenceless terms.  On the
other hand, a calculation shows that, for an arc of a trajectory of
$\vec{P_1}$ with $r > 2$ and $z \in [-1/3,1/3]$, $\frac{d\theta}{dz}$ is
$3\pi b(\frac{1+3z}{2})$
greater than it is for a similar arc in $\vec{P_2}$, and the integral over
$z$ of this difference is $2\pi$. In other words, a leaf of
${\cal P}_1$ with $r>2$ has the same endpoints as some leaf of ${\cal P}_2$, but
winds in the $\theta$ direction by an extra $2\pi$ exactly.  If we
concatenate the mirror image of ${\cal P}_2$ to ${\cal P}_1$ in the manner of the
mirror-image construction, the result is that leaves with $r>2$ wind an angle
of $2\pi$ while leaves with $r<2$ wind an angle of $0$. Therefore for the
plug $\cal P$, the two sides of the circle $\alpha(c)$ do not wind around the
same amount, and $\alpha(c)$ is not a meridian for $\cal P$, as shown in
Figure~\ref{fmeridian}b.  In fact, $\cal P$ is integrally Dehn-twisted.

The plug $\cal P$ has two annuli of closed leaves.  Since the stopped set of
$\cal W$ is a circle, $\cal W$ can be inserted in such a way that all of these
closed leaves are broken.   The insertion of $\cal W$ into $\cal P$ produces
the desired plug $\cal D$ with two closed leaves, as desired.
\end{proof}

To achieve a measured foliation with very few closed leaves, namely two, on an
arbitrary compact 3-manifold $M$, we can insert a single copy of $\cal W$
that breaks all closed leaves of the foliation of Theorem~\ref{thdiscrete}.
As an alternative to the proof of Lemma~\ref{ldehn}, we could equally well
insert copies of the plug $\cal P$ to effect Dehn surgery on $T^3$ or $\R P^2
\times S^1$ and then use one copy of $\cal W$ in the final step.

\subsection{The PL case}

The foliation $\cal T$ on $T^3$ is also a measured PL foliation. An irrational
rotation of $\R P^2$ can also be realized as an area-preserving PL
homeomorphism.  Therefore the following lemma establishes
Theorem~\ref{thdiscrete} by the same reasoning as in the smooth case:

\begin{lemma} 
There exists a PL, measured, integrally Dehn-twisted plug $\cal D$
with two closed leaves.
\end{lemma}
\begin{proof}
Let $R$ be a compact, PL submanifold of the plane, given in coordinates $x$ and
$y$.  Let $f:R \to R$ be a PL homeomorphism, and let $l$ be a real number.  Let
$\cal L$ be the foliation of $R \times I$ such that, for fixed $(x,y) \in R$ and $z
\in \R$, the set $\{(x,y + lx,z) | (x,y+lz) \in R, z \in I\}$ is a leaf.  Orient
the leaves in the direction of increasing $z$.  The {\em slanted suspension} of
$f$ with slope $l$ is defined as the space $R \times I$ with $(x,y,1)$ identified
to $(f(x,y),1)$, together with the foliation $\cal S$ induced by $\cal L$.  The
slanted suspension $\cal S$ is manifestly a PL foliation also.  Moreover, if $f$
preserves area, then $\cal S$ is measured.

Let $T$ be the trapezoid in the plane with vertices $a_1 = (0,0)$, $a_2 =
(0,2)$, $a_4 = (1,1)$, and $a_5 = (1,0)$, and let $a_3 = (\frac13,0)$, as shown
in Figure~\ref{ftrap}.  Let $U$ be the reflection of $T$ about the line $x =
\frac12$, and let $b_i$ be the image of $a_i$ under this reflection. Let $f:T
\to U$ be the unique PL map which sends $a_i$ to $b_{6-i}$ and is linear on each
of the three triangles which share the vertex $a_3$.  Evidently $f$ is an
area-preserving PL homeomorphism.  Moreover, $f$ decreases the $y$ coordinate of
$a_2$ the most and increases the $y$ coordinate of $a_4$ the most.

\fullfigure{}{The map $f:T \to U$}
 
Let $R_1 = [-1,1] \times [0,3]$ be a rectangle consisting of four congruent
copies of $T$ as shown in Figure~\ref{fg1}.  We can conjugate $f$ with three
isometries of the plane to extend $f$ to an area-preserving PL homeomorphism
$g_1:R_1 \to R_1$ as also indicated in Figure~\ref{fg1}. It is easy to
check that the slanted suspension ${\cal S}_1$ of $g_1$ with slope $1$ is a PL
analogue of the semi-plug ${\cal W}_s$; it is a semi-plug
that stops a circle and has one closed leaf.

\fullfigure{}{The map $g_1$}

Let $R_2$ be the rectangle $R_1$ union a triangle with vertices $(0,3)$,
$(1,3)$, and $(1,5)$.  Let $g_2:R_2 \to R_2$ be a map pieced together from
four copies of $f$:  $f$ itself, its reflection (in the sense of conjugation)
about the $y$ axis, the rotation by 180 degrees of that reflection about the
point $(-\frac12,\frac32)$, and the image of $f$ under the affine
transformation $(x,y) \mapsto (1-x,5-2x-y)$.  The four copies of $f$
determine $g_2$ everywhere except in the triangle with vertices $(0,2)$,
$(3,1)$, and $(3,3)$; $g_2$ is defined by $g_2(x,y) = (x,y+x)$ on this
triangle.  Figure~\ref{fg2} gives a diagram of the map $g_2$.  Let ${\cal P}_2$
be the slanted suspension of $\alpha_2$ with slope $1$.  It is easy to check
that ${\cal P}_2$ is a semi-plug with one closed leaf as well.

\fullfigure{}{The map $g_2$}

The goal is to concatenate ${\cal S}_1$ and ${\cal S}_2$ to produce a plug ${\cal
S}_{PL}$, but we must be careful to properly match the entry and exit regions. 
Let $F_{1,\pm}$ and $F_{2,\pm}$ be the entry and exit regions of ${\cal S}_1$ and
${\cal S}_2$; all four are subsets of $\R^2 \times [0,1]$ with $\R^2 \times \{0\}$
and $\R^2 \times \{1\}$ identified.  Let $F = F_{1,-} = F_{2,-}$, and let
$\pi_1:F_{1,+} \to F$ and $\pi_2:F_{2,+} \to F$ be vertical projections.  Note
that $\pi_1$ and $\pi_2$ preserve transverse measure and that transverse measure
agrees on $F_{1,-}$ and $F_{2,-}$.  The essential property of ${\cal S}_1$ and
${\cal S}_2$ is that if $p$ and $\pi_1^{-1}(q)$ are endpoints of a leaf of ${\cal
S}_1$, then $p$ and $\pi_2^{-1}(q)$ are the endpoints of another leaf of ${\cal
S}_2$ that winds around the suspension direction one extra time.  Therefore if
${\cal S}_1$ and $\bar{{\cal S}_2}$ are concatenated with the map $\pi_2 \circ
\pi_1^{-1}$ as a gluing map, the result is an integrally Dehn-twisted plug $\cal
S$ whose geometry is very similar to that of the plug $\cal P$ constructed in the
previous subsection.  The plug $\cal S$ has two closed leaves as desired.
\end{proof}

\subsection{Non-compact 3-manifolds}
\label{snoncompact}

The only extra difficulty in establishing Theorem~\ref{thdiscrete} in the
non-compact case is that the Lickorish-Wallace theorem does not generalize. It
is not true, for example, that every open 3-manifold is obtained from $\R^3$ by
Dehn surgery on a locally finite link, because any such surgery produces a
manifold with only one end. (Recall that the set of ends of a manifold is the
inverse limit of the connected components of complements of compact subsets.)
On the other hand, the theorem does have the following useful generalization

\begin{theorem} (Generalized Lickorish-Wallace theorem) Given two 3-manifolds
$A$ and $B$ with a homeomorphism $\alpha:\d A \to \d B$ there exists
a 3-manifold $\hat{A}$ obtained from $A$ by integral Dehn surgery on a
link disjoint from $\d A$ such that $\alpha$ extends to a homeomorphism
$\alpha:\hat{A} \to B$. \label{thglw}
\end{theorem}
\begin{proof} (Sketch) For closed manifolds, the Lickorish-Wallace theorem
essentially says that any 3-manifold bounds a 4-manifold, since an integral
Dehn surgery is a Morse reconstruction at a critical point of index 2, for
3-manifolds viewed as level surfaces of Morse functions on 4-manifolds. 
Contrariwise, if two 3-manifolds are level surfaces of the same Morse
function, the Morse stratification produces a Dehn surgery connecting them,
since any possible Morse reconstruction can be reproduced with Dehn surgery. 
In the case at hand, $C = A \cup B \cup (\d A \times I)$ is a closed
3-manifold if $\d A \times \{0\}$ is identified with $\d A$ and $\d A \times
\{1\}$ is identified with $\d B$ using $\alpha$. The manifold $C$ bounds a
4-manifold $W$; in fact, $W$ is naturally an orthant manifold if $\d A \times
I$ is positioned to meet $A$ and $B$ orthogonally.  Choosing a Morse function
that is 0 on $A$, 1 on $B$, and $x$ on $\d A \times \{x\}$, we obtain a
sequence of Morse moves that connect $A$ to $B$, which can again be converted
to Dehn surgeries.
\end{proof}

\begin{lemma} Every non-compact, orientable 3-manifold $M$ can be realied by a
Dehn surgery on a locally finite link in an infinite collection of spheres with
connecting handles.  \label{lopensur}
\end{lemma}
\begin{proof}

Let $M$ be a non-compact 3-manifold.  Consider a locally finite collection of
embedded surfaces in $M$ which separate $M$ into compact 3-manifolds
$M_1,M_2,\ldots$ with boundary. Let $S_1,S_2,\ldots$ be a collection of
3-spheres, and connect $S_i$ to $S_j$ by a handle (in the sense of connected
sums) for each connected component of $M_i \cap M_j$.  The resulting manifold
$P$ is tiled by punctured 3-spheres $P_1,P_2,\ldots$ with the property that
$P_i \cap P_j$ consists of $n$ 2-spheres if $M_i \cap M_j$ has $n$ connected
components.  By attaching handles to these 2-spheres, we can obtain a new
tiling $P'_1,P'_2,\ldots$ such that $P'_i \cap P'_j$ is homeomorphic to $M_i
\cap M_j$.  Applying Theorem~\ref{thglw} to $P$, there exists a finite
surgery in each $P'_i$ that yields $M_i$.  The union of all such surgeries is
a locally finite surgery on $P$ that yields $M$.
\end{proof}

Given Lemma~\ref{lopensur}, Theorem~\ref{thdiscrete} is established with the aid
of a volume-preserving, twisted plug $\cal H$ with base $D^2 \amalg D^2$ and
with support $(D^2 \times I) \# (D^2 \times I)$.  In particular, the base of
$\cal H$ is not connected, but the support is; $\cal H$ is therefore a {\em
wormhole plug}.  The insertion of $\cal H$ into a foliation of a disconnected
manifold effects a connected sum between two different components of the
manifold.

\begin{lemma}  There exists a smooth, measured plug $\cal H$ with
base $D^2 \amalg D^2$ and support $(D^2 \times I) \# (D^2 \times I)$.
\end{lemma}
\begin{proof}
The first step in constructing $\cal H$ is to construct a semi-plug ${\cal
H}_s$ with the same support.  Figure~\ref{fworm} shows a flow which is a
2-dimensional analogue of the desired semi-plug:  It is a flow in an orthant
manifold which is homeomorphic to an annulus.  The outside boundary is a
square, the inside boundary is an inverted square, and both the annulus and its
flow might be invariant under inversion in a circle in the center of the
annulus.  Excepting the two fixed points, it is otherwise a flow bordism. 
Roughly speaking, the 3-dimensional semi-plug ${\cal H}_s$ is parallel to a
flow obtained by revolving the 2-dimensional analogue about the vertical axis
and adding motion in the angular direction to turn the fixed points into a
closed trajectory.

\fullfigure{}{A 2-dimensional wormhole}

Explicitly, parameterize $\R^3$ by Cartesian coordinates $x$, $y$, and $z$,
and let $C$ be the cylinder given by $x^2 + y^2, z^2 \le 16$.  Let
$\omega = dx \wedge dy$ be a flux form on $C$; the parallel foliation
is simply parallel vertical line segments.  Let $\alpha:\R^3 \to \R^3$
be given by
$$\alpha(x,y,z) = {1 + \frac{1}{x^2 + y^2 + z^2}}(x,y,z).$$
Define ${\cal H}_s$ to be the foliation parallel to the flux form
$$\alpha^*(\omega) + x dx \wedge dz + y dy \wedge dz$$
on the domain $\alpha^{-1}(C)$.  The foliation ${\cal H}_s$ has
all of the claimed properties.

The mirror-image construction applied to ${\cal H}_s$ produces a plug ${\cal H}_m$
whose support $H_m$ consists of two cubes connected by two handles, rather than
the desired two cubes connected by one handle.  However, the manifold $H_m$ can
be written as 
$$H_m \cong (D^2 \times I) \# (D^2 \times I) \# (S^2 \times S^1).$$
Since there exists a Dehn surgery on $S^2 \times S^1$ that yields $S^3$, there
is a way to insert copies of the Dehn-twisted plug $\cal D$ into ${\cal H}_m$ to
produce a plug with support $(D^2 \times I) \# (D^2 \times I)$.  This plug is
$\cal H$.
\end{proof}

A PL analogue of $\cal H$ also exists; the details are omitted.

\section{No closed trajectories}

Since the construction of the proof theorem~\ref{thmain} is a modification 
of a Schweitzer plug, we begin with a brief review of that example.

\subsection{Schweitzer's construction}

\label{sschweitzer}

If $a$ and $b\ne 0$ are real numbers, let $a \bmod b$ be the
corresponding element in the circle $\R/b\Z$.  Let $\tau$ be an irrational real
number. Let 
$$w(x) = \frac1\pi\left(\tan^{-1}(x+1) - \tan^{-1}(x)\right)$$
and consider a sequence of open intervals $I_n$ of length $|I_n| = w(n)$ placed
on the unit circle in the same order as $n \bmod \tau$.  More specifically, let
$$I_n = (a_n \bmod 1, a_n+w(n) \bmod 1) \subset \R/\Z = S^1,$$ where
$$a_n = \sum_{k:k \bmod \tau \in [0,n \bmod \tau)} w(k).$$
Since the total length of the $I_n$'s is exactly 1, they are dense in $S^1$
but they do not intersect each other.  Because of the ordering of the 
$I_n$'s, there exists a homeomorphism $\alpha:S^1 \to S^1$ such that
$\alpha(I_n) = I_{n+1}$.  The map $\alpha$ is a {\em Denjoy homeomorphism}.
It is realized as a $C^1$ diffeomorphism if its derivative of
on $I_n$ is defined to be
$$\frac{d \alpha}{dx} = 1 + \frac{|I_{n+1}| - |I_n|}{|I_n|} b(L_n(x)),$$
where $L_n:I_n \to [0,1]$ is a linear isomorphism, for those $n$ such that
$4|I_{n+1}| > 3|I_n|$.  For the finite number of $n$ such that this inequality
fails, let $\alpha$ be any diffeomorphism from $I_n$ to $I_{n+1}$ with
derivative 1 at the endpoints. The derivative of $\alpha$ is 1 outside of the
$I_n$'s, and since
$$\lim_{n \to \pm \infty} \frac{|I_{n+1}| - |I_n|}{|I_n|} = 0,$$
the derivative is continuous.  The map $\alpha$ has no periodic
orbits and has a unique minimal set, namely $S^1 - \bigcup_n I_n$.

Let $\cal D$ be the suspension of $\alpha$ and let $T$ be its support.  The
manifold $T$, which is a torus because $\alpha$ preserves
orientation, is {\em a priori} only a $C^1$ manifold, but it has a smooth
refinement such that $\cal D$ is parallel to a $C^1$ vector field $\vec{D}$.
Let $m$ be the subset of $T$ which is the suspension of the minimal
set $S^1 - \bigcup_n I_n$; $m$ is the minimal set of $\cal D$.  The set
$m$ is a {\em Denjoy continuum}.

Consider the manifold $T \times [-1,1]$ with the $[-1,1]$ factor parameterized
by $z$.  Let $f$ be a non-negative, non-zero $C^1$ function on $T$ which
vanishes on $m$.  Consider a vector field $\vec{E}$ given by the formula:
$$\vec{E} = \vec{D} + z^2\frac{\d}{\d z} + f\frac{\d}{\d z}.$$ 
By inspection of $\vec{E}$, the parallel foliation $\cal E$ has no closed leaves,
because all leaves either travel in the positive $z$ direction or coincide with
leaves of $\cal D$ on $T \times \{0\}$.  Moreover, $\cal E$ has infinite
leaves contained in its minimal set $m \times \{0\}$.  Therefore, by
lemma~\ref{linfleaf}, $\cal E$ is a semi-plug whose base is a torus.  The
mirror-image construction applied to $\cal E$ yields a plug $\cal F$ which also
has no closed leaves and has the same base.  We identify the support of $\cal F$ with
$T$.

Since the torus has no boundary, $\cal F$ is not insertible.  However, $\cal F$
also possesses a leaf $l$ with two endpoints; we can take $l$ to be the
extension of the leaf in $\cal E$ containing $(p,0)$, where $p \in T$ satisfies
$f(p) > 0$.  The leaf $l$ has a foliated tubular neighborhood $N_l$ consisting
entirely of leaves with two endpoints.  The restriction $\cal S$ of $\cal F$ to
$T \times I - N_l$ is therefore a plug with base $pT$, a punctured
torus.  Since $l$ is unknotted (which follows from the fact that $\vec{E}$ is
non-negative in the $z$ direction, and,  by the mirror image construction,
the leaves in $N_l$ do not twist around $l$),  $\cal S$ is an
untwisted plug.  Following Section~\ref{sintro}, copies of $\cal S$ can be
inserted to break any discrete collection of closed leaves in a foliation. 
Indeed, as discussed in the appendix, a 3-dimensional plug with a twisted or
knotted leaf neighborhood removed can nevertheless be extended to an untwisted,
insertible plug. The existence of the plug $\cal S$ together with Wilson's
theorem (or its variant in Section~\ref{sdiscrete}) establishes a counterexample
to the usual Seifert conjecture for all 3-manifolds \cite{Schweitzer}.

\subsection{Preserving volume}

The difficulty in making Schweitzer's construction volume-preserving is the
fact that $\cal D$ does not possess a transverse measure in the sense of
Section~\ref{sintro}. Such a measure would induce an $\alpha$-invariant measure
$\mu$ on the circle which is locally equivalent to Lebesgue measure. By
compactness, the total $\mu$-measure of the circle would be finite, but the
$I_n$'s would have equal and non-zero measure, a contradiction.  In other
words, any homeomorphism of $S^1$ conjugate to $\alpha$ has the inevitable
effect of squeezing $I_n$ as $n$ goes to $\infty$ and stretching $I_n$ as $n$
comes from $-\infty$.  Our strategy for overcoming this difficulty is to
compensate squeezing of $I_n$ by stretching in the $z$ direction.  This
transverse stretching must be sufficiently slight that there is no net motion
in the negative $z$ direction.

Finding a suitable amount of transverse stretching is the difficult part of
Theorem~\ref{thmain} because it is bounded both above and below by different
constraints in the construction.  One particular problem is that, if the
rotation number $\tau$ of the Denjoy homeomorphism $\alpha$ is approximated too
closely by rationals, orbits of rotation by $\tau$ are too unevenly distributed
for the construction to work.  Although any irrational number whose continued
fraction expansion has bounded coefficients would work in principle, we let
$\tau = \frac{1+\sqrt{5}}{2}$ be the golden ratio for simplicity. In any case,
the construction requires some involved if elementary Diophantine estimates. 
For convenience, the presence of the constant $C$ in an equation will mean that
there exists a real number $C > 0$ such that the equation holds.  Although $C$
is independent of all variables, it may have a different value in different
equations, or even in different sides of the same equation.

The first step is to construct the Denjoy foliation, or at least its minimal
set, in such a way that the underlying torus has a convenient measure. 
Consider the cylinder $S^1 \times \R$ parametrized by $\theta$ and $\phi$. Let
$S_n \subset S^1 \times \R$ be a sequence of infinite cylindrical strips such
that the intersection $S_n \cap (S^1 \times \{\phi\})$ has length $w(n-\phi)$
and such that the $S_n$'s have the same ordering.  Explicitly, let $I_{n,\phi}
\subset S^1 \times \{\phi\}$ be the open interval given by
\begin{multline*}
I_{n,\phi} = (a_{n,\phi} \bmod 1, a_{n,\phi}+w(n-\phi) \bmod 1) \\
    \subset S^1\times \R = \R/\Z\times \R,
\end{multline*}
where
$$a_{n,\phi} = \sum_{k \bmod \tau \in [0,n \bmod \tau)} w(k-\phi),$$
and let $S_n$ be the union of all intervals $I_{n,x}$.  
Let $\sigma:S^1 \times \R
\to S^1 \times \R$ be the homeomorphism given by 
$$\sigma(\theta,\phi) = (\theta+a_{1,\phi+1},\phi+1).$$
The map $\sigma$ is smooth, preserves area on $S^1 \times \R$, and sends $S_n$
to $S_{n+1}$.  The quotient $T$ of $S^1 \times \R$ by $\sigma$ has an open
strip $S$ which is the image of each $S_n$ under the quotient map, and the
volume form $d\theta \wedge d \phi$ descends to a form $\mu$ on $T$.  The
complement $m$ of $S$ in $T$ is clearly a Denjoy continuum.  In fact, by this
definition, the pair $(T,m)$ is an explicit smooth refinement of the objects of
subsection \ref{sschweitzer} with the same name.

Consider $S^1 \times \R \times [-1,1]$ with the third coordinate parameterized
by $z$ and with measure $\mu \wedge dz$.  The next step is to define vector fields
$\vec{h}$ and $\vec{v}$
on $S^1 \times \R \times [-1,1]$ with the following properties:
\begin{itemize}
\item[(i)] They are both invariant under $\sigma \times id$ and therefore descend to $T
\times [-1,1]$.
\item[(ii)] They are both divergenceless $C^1$ vector fields (relative to
$d\theta \wedge d\phi \wedge dz$ or $\mu \wedge dz$) whose $\phi$
components vanish.
\item[(iii)] The vector field $\vec{h}+\frac{\d}{\d \phi}$ is parallel to $m
\times \{0\}$.  On the other hand, $\vec{v}$ vanishes on $m \times \{0\}$.
\item[(iv)] The $z$ component of $\vec{v}$ is positive except on
$m \times \{0\}$ and exceeds the absolute value of the $z$ component
of $\vec{h}$.
\end{itemize}
Assuming for the moment the existence of $\vec{v}$ and $\vec{h}$,
the trajectories of the vector field $\vec{E}' = \vec{v}+\vec{h}+\frac{\d}{\d\phi}$
have the same geometry as those of $\vec{E}$; moreover $\vec{E}'$ is
divergenceless.  Following the rest of Schweitzer's construction and the
formalism of Section~\ref{splugs}, the flux form given by $\vec{E}'$ yields a
measured $C^1$ plug, which establishes Theorem~\ref{thmain}.

We temporarily fix a value of $\phi$ and work in the coordinates $\theta$ and
$z$ with the measure $d\theta \wedge dz$. Let $w'(x)$ denote the derivative
of $w(x)$.  Let
$$f(\theta) = \frac{w'(n-\phi)}{w(n-\phi)} b(L_{n,\phi}(\theta))$$
for $\theta \in I_{n,\phi}$ and $0$ elsewhere, where $L_{n,\phi}:I_{n,\phi}
\to [0,1]$ is a direction preserving linear isomorphism.  Let
$$F(\theta) = w(n-\phi)^{3/2} B(L_{n,\phi}(\theta))$$
for $\theta \in I_{n,\phi}$ and $0$ elsewhere.  Define $\vec{h}$ and $\vec{v}$
by the equations
\begin{eqnarray}
H(\theta,z) & = &\frac12\int_{\theta-z}^{\theta+z}
\int_0^{\theta_1} f(\theta_2) d\theta_2 d\theta_1 \label{ebighdef}\\
V(\theta,z) & = & \frac C z \int_{\theta-5z}^{\theta+5z}
\int_0^{\theta_1} F(\theta_2) d\theta_2 d\theta_1 \label{ebigvdef}\\
\vec{h} &=& J(\vec{\nabla} H)=
(-\frac {\partial H} {\partial z},\frac {\partial H} {\partial\theta} ) \nonumber \\
\vec{v} &=& J(\vec{\nabla} V)=
(-\frac {\partial V} {\partial z},\frac {\partial V} {\partial\theta} ) \nonumber
\end{eqnarray}
extended to $z = 0$ by continuity.

Except for $C^1$ continuity, properties (i), (ii), and (iii) are routine.
The fact that $\vec{h}$ is $C^1$ follows from $C^2$ continuity of $H$,
which is immediate from the continuity of $f$.  The function $F$
is $C^1$ as follows:  The derivative exists on each $I_{n,\phi}$,
and it extends continuously to a function $\tilde{F}:S_1 \to \R$
which is zero outside of the $I_{n,\phi}$'s (check).  We claim
that $F$ is the antiderivative of $\tilde{F}$.  It could only
disagree with the antiderivative if it were discontinuous or
if the set $F(S_1 - \bigcup_n I_{n,\phi})$ had
non-zero Lebesgue measure, and neither of these is the case.
Since $F$ is $C^1$, $V$ is $C^3$ everywhere except where $z = 0$; at such
points $V$ is $C^2$ by L'Hospital's rule.  Therefore $\vec{v}$ is $C^1$
also.

Property (iv) is the heart of the matter, and we prove it with a sequence of
lemmas.  If $a,b \in \R/\tau\Z$, let $Z(a,b)$ be the set of all integers $n$
such that $n \bmod \tau \in (a,b)$, and let $d(a,b)$ be the distance from $a$
to $b$ on the circle $\R/\tau\Z$.  Here the notation $(a,b)$ denotes an
interval $(a,b)$ whose endpoints are $a$ and $b$ and which is oriented from
$a$ to $b$ in the natural orientation of the circle.

\begin{lemma} Let $F_n$ be the $n$th Fibonacci number, with
$F_0 = F_1 = 1$ and $F_{n+2} = F_{n+1} + F_n$.  Then
$$d(F_n \bmod \tau, 0) = \tau^{-n}.$$
Moreover, $F_{2n} \bmod \tau $ and $F_{2n+1} \bmod \tau$ converge to 0 from
opposite sides.  \label{lfibhowgood}
\end{lemma}
\begin{proof}
By induction and applying the identity $\tau^{-1} = \tau - 1$, we have
\begin{align*}
F_{2n} \bmod \tau &= (F_{2n-1}+ F_{2n-2})\bmod \tau \\
    &= -\tau^{-(2n-1)} + \tau^{-(2n-2)} \\
    &= \tau^{-(2n-1)}(\tau -1) = \tau^{-2n},
\end{align*}
and
\begin{align*}
\tau - F_{2n+1} \bmod \tau &= \tau - (F_{2n}+ F_{2n-1})\bmod \tau \\
    &= \tau  -(\tau^{-2n} - \tau^{-(2n-1)} +\tau ) \\
    &= \tau^{-(2n-1)} - \tau^{-2n} \\
    &= \tau^{-2n}(\tau - 1)=\tau^{-(2n+1)}.
\end{align*}
\end{proof}
\begin{lemma} If $0 < p < F_n$, then
$$d(F_n \bmod \tau,0) < d(p \bmod \tau,0).$$ \label{lfibopt}
\end{lemma}
\begin{proof} Applying the identity $\tau^{-1} = \tau - 1$ inductively to 
$\tau^{-n}$ yields
$$\tau^{-n} = (-1)^n (F_n - F_{n-1} \tau).$$
Hence, the lemma can be rephrased as
$$|F_n - \tau F_{n-1}| < |p - \tau q|,$$
for some integer $q$. The proof follows now from Theorem 182 of Hardy and Wright~\cite{Hardy}, since the
ratios $\frac{F_{n-1}}{F_n}$ are the partial evaluations of the continued
fraction expansion of $\tau^{-1}$.
\end{proof}

\begin{lemma} Let $n_1$ and $n_1+k_1$ be a pair of consecutive elements in
$Z(a,b)$, and let $n_2$ and $n_2+k_2$ be another such pair.  Then
$k_1 < Ck_2$. \label{lquasilat}
\end{lemma}
\begin{proof}
The case in which $(a,b)$ is more than half of the circle $\R/\tau\Z$ is
trivial.  In the non-trivial case, $d(a,b)$ is the length of the interval
$(a,b)$.  Choose the largest $n$ such that $d(F_n \bmod \tau,0) >
d(a,b)$.  Since $k_1 \bmod \tau < d(a,b)$, it follows that
$k_1 > F_n$ by Lemma~\ref{lfibopt}.  On the other hand, by
Lemma~\ref{lfibhowgood}, since $d(F_{n+1} \bmod \tau,0) < d(a,b)$, $d(F_{n+3}
\bmod \tau,0) < d(a,b)/2$ and $d(F_{n+4} \bmod \tau,0) < d(a,b)/2$ also, and
$F_{n+3} \bmod \tau$ and $F_{n+4} \bmod \tau$ are on opposite sides of 0.
At least one of $k_1+F_{n+3} \bmod \tau$ and $k_1+F_{n+4} \bmod \tau$ is in 
$(a,b)$, since $k_1 \bmod \tau$ is at least $d(a,b)/2$ away from
one endpoint of $(a,b)$.  Therefore $F_n < k_1 \le F_{n+4}$,
and since $F_{n+4} < C F_n$ and all arguments also apply to $k_2$,
the conclusion follows.
\end{proof}

\begin{lemma}  If $0 \notin Z(a,b)$, then
$$\sum_{n \in Z(a,b)} \frac{1}{|n|^3} \le 
C\frac{\displaystyle \sum_{n \in Z(a,b)} \frac{1}{|n|^5}}
    {\displaystyle \sum_{n \in Z(a,b)} \frac{1}{n^2}}.$$
\label{l352}
\end{lemma}
\begin{proof}
Let $k$ be the element of $Z(a,b)$ with the least absolute value, and
assume without loss of generality that $k > 0$.  By Lemma~\ref{lquasilat},
the minimum gap between elements of $Z(a,b)$ is at least $Ck$.
It follows that
$$\sum_{n \in Z(a,b)} \frac{1}{n^2} \le 
    2\sum_{n=0}^\infty \frac{1}{(k + nCk)^2} = 
    \frac{2}{\displaystyle k^2 \sum_{n=0}^\infty \frac{C^2}{(C+n)^2}}
    = \frac{C}{k^2}.$$
Therefore,
$$\frac{\displaystyle \sum_{n \in Z(a,b)} \frac{1}{|n|^5}}
    {\displaystyle \sum_{n \in Z(a,b)} \frac{1}{n^2}} > C\sum_{n \in Z(a,b)}
    \frac{k^2}{|n|^5} > C\sum_{n \in Z(a,b)}\frac{1}{|n|^3}.$$
\end{proof}

\begin{corollary} For any $a$, $b$, and $\phi$,
$$\sum_{n \in Z(a,b)} w(n-\phi)^{3/2} \le
    C\frac{\displaystyle \sum_{n \in Z(a,b)} w(n-\phi)^{5/2}}
    {\displaystyle \sum_{n \in Z(a,b)} w(n-\phi)}.$$
\label{c352}
\end{corollary}
\begin{proof}
The case $0 \in Z(a,b)$ is trivial; suppose that $0 \notin Z(a,b)$.
Without loss of generality, $0 \le \phi < 1$, and with this restriction,
$$\frac{C}{n^2} > w(n-\phi) > \frac{C}{n^2}.$$
The inequality renders the corollary equivalent to Lemma~\ref{l352}.
\end{proof}

\fullfigure{}{Bounding the length of $I_v'$}

\begin{lemma} If $n_1$ and $n_2$ are distinct integers, there exists 
an integer $n_3$ such that $n_3 \bmod \tau \in (n_1 \bmod \tau, n_2 \bmod \tau)$
and
$$|n_3| < C \max(|n_1|,|n_2|)$$
\label{l312}
\end{lemma}
Lemma~\ref{l312} is a corollary of Lemmas~\ref{lfibhowgood} and \ref{lfibopt} in
the same way as Lemma~\ref{lquasilat} is.

\widefullfigure{}{Isotopy of a plug insertion}

\begin{corollary} If $n_1$ and $n_2$ are distinct integers, then the
intervals $I_{n_1,\phi}$ and $I_{n_2,\phi}$ satisfy
$$d(I_{n_1,\phi},I_{n_2,\phi}) > C\max(|I_{n_1,\phi}|,|I_{n_2,\phi}|).$$
\label{c312}
\end{corollary}
\begin{proof}
Combining the estimate
$$\frac{C}{n^2} > w(n - \phi) > \frac{C}{n^2}$$
with Lemma~\ref{l312}, there exists an $n_3$ such that $I_{n_3,\phi}$
lies between $I_{n_1,\phi}$ and $I_{n_2,\phi}$ and such that
$$w(n_3 - \phi) > C \min(w(n_1 -\phi),w(n_2 - \phi)).$$
Since the length of $I_{n,\phi}$ is $w(n-\phi)$, the lemma follows
from the fact that $I_{n_1,\phi}$ and $I_{n_2,\phi}$ are sufficiently
far apart to make room for $I_{n_3,\phi}$.
\end{proof}

To establish property (iv), expand $\vec{v}$ and $\vec{h}$ as
\begin{eqnarray*}
\vec{v} & = & v_\theta \frac{\d}{\d \theta} + v_z \frac{\d}{\d z} \\
\vec{h} & = & h_\theta \frac{\d}{\d \theta} + h_z \frac{\d}{\d z}.
\end{eqnarray*}
It suffices to show that $v_z > |h_z|$ when $z \ne 0$.  Once again fix $\phi$,
and note from equations~\ref{ebighdef} and~\ref{ebigvdef} that these two
quantities are given by
\begin{eqnarray*}
h_z & = & \frac12\int_{\theta-z}^{\theta+z} f(\theta_1) d\theta_1 \\
v_z & = & \frac C z \int_{\theta-5z}^{\theta+5z} F(\theta_1) d\theta_1.
\end{eqnarray*}
The absolute value $|h_z|$ is also bounded by
$$h_{\abs} = \frac12\int_{\theta-z}^{\theta+z} |f(\theta_1)| d\theta_1. $$
Let $I_v = [\theta-5z,\theta+5z]$ be the domain of integration for $v_z$,
and similarly let $I_h = [\theta-z,\theta+z]$.  One possibility is that
$I_h$ is a subset of some $I_{n,\phi}$.  In this case, the inequalities
$$z < |I_{n,\phi}|$$
$$B(x) > b(x)$$
$$Cw(n-\phi)^{3/2} > |w'(n-\phi)|$$
together imply that $v_z > h_{\abs}$.  

Alternatively, suppose that $I_h$ does not lie in a single $I_{n,\phi}$.  Let
$I_v'$ be the closure of the union of all $I_{n,\phi}$'s which are
contained in $I_v$.  Since $I_h$ is the middle fifth of $I_v$, and since the
integrand of $h_{\abs}$ is only non-zero in the middle third of an interval
$I_{n,\phi}$, the integrand of $h_{\abs}$ is zero in the region in $I_h -
I_v'$.  \Ie,
\begin{equation}
h_{\abs} \le \int_{I_v'} |f(\theta_1)| d\theta_1. \label{ehbound}
\end{equation}
The region $I_v - I_v'$ in general consists of a subinterval of some
interval $I_{n_1,\phi}$ on one side and a subinterval of some other interval
$I_{n_2,\phi}$ on the other side.  By hypothesis, $I_v'$ contains at least
one point in $I_h$, which is the middle fifth of $I_v$, and therefore if
$|I_v'| < C|I_v| = Cz$, then the intervals $I_{n_1,\phi}$ and
$I_{n_2,\phi}$ both have length at least $Cz$.  (See Figure~\ref{fintervals}.)
It follows by Corollary~\ref{c312} that $|I_v'| > Cz$ and that
\begin{equation}
v_z > \frac C {|I_v'|} \int_{I_v'} F(\theta_1) d\theta_1. \label{evbound}
\end{equation}

In the main case, the set $Z_v$ of all $n$ such that $I_{n,\phi} \subset I_v'$
is exactly a $Z(a,b)$.   In this case, the right sides of 
equations~\ref{ehbound} and~\ref{evbound} are related by Corollary~\ref{c352},
which demonstrates that $v_z > h_{\abs}$, as desired.  The alternative
possibilities are that $a = n_a \bmod \tau$ or that $b = n_b \bmod \tau$ and that
$Z_v$ is $Z(a,b)$ union $\{n_a\}$ or $\{n_b\}$ or both; these exceptional cases
can be treated in the same way as the main case.

\fullfigure{}{A plug with a knotted hole}

\section{Appendix:  Plugs with knotted holes}

The Schweitzer construction and its modifications here and in the work of
Harrison~\cite{Harrison} suggest but do not depend on the following question:
Suppose that $\cal P$ is a plug whose base is a closed, oriented surface and
suppose that $\cal P$ has a knotted leaf with two endpoints.  Let $N_l$ be a
foliated tubular neighborhood of $l$ and let ${\cal P}_l$ be $\cal P$ with
$N_l$ removed.  Since ${\cal P}_l$ is twisted, can it be inserted into a
foliation of $M$ without changing the topology of $M$?

In Harrison's construction, the base of $\cal P$ is a torus and, moreover,
$\cal P$ is a slanted suspension (in the $C^2$ category) of a homeomorphism of
an annulus, made into a plug with the mirror-image construction.  In this case,
following Harrison, there exists a finite cover of $\cal P$ such that a lift of
$l$ is necessarily unknotted, and by the mirror-image construction $N_l$ is
necessarily untwisted.

Taking the general case, suppose that the base of $\cal P$ is $S$; the base of
${\cal P}_l$ is therefore $pS$, or $S$ with one puncture.  Let $P_l$ be the
support of ${\cal P}_l$.  Consider Schweitzer's insertion of ${\cal P}_l$ into
the un-plug with base a disk $D^2$ in Figure~\ref{fstandard}a.  If ${\cal P}_l$
were an untwisted plug, its insertion would be realized by an embedding
$\alpha$ of $pS \times I$ which is a thickening of the insertion map for the
base.  Ignoring the vertical foliation on $D^2 \times I$, this embedding is
isotopic to the standard embedding of the closed surface $S$, punctured and
thickened, as shown in Figure~\ref{fstandard}b and Figure~\ref{fstandard}c. 
The complement $(D^2 \times I) - \alpha(pS \times I)$ is topologically the
exterior of a solid torus connected by a handle to a solid torus. 
Recognizing $P_l$ as $S \times I$ with a knotted hole, it admits an embedding
$\beta$ in $D^2 \times I$ whose complement is homeomorphic to the
compement of $\alpha(pS \times I)$, as shown in Figure~\ref{fknotted}.
Since the complements are the same, the insertion of ${\cal P}_l$ does
not change the topology of $D^2 \times I$ even thought ${\cal P}_l$ is
twisted.

\end{document}